\documentclass[11pt,leqno,a4paper]{article}

\usepackage[hypertex]{hyperref}
\usepackage{amsfonts}
\usepackage{latexsym}
\usepackage{amssymb,amsmath}
\usepackage{epsfig}

\def\intP{\stackrel{\circ}{P}}
\def\va{\hfill\break}

\def\to{\rightarrow}

\def\ep{\varepsilon}
\def\H{\mathbb H}
\newcommand{\C}{\mathbb C}
\newcommand{\R}{\mathbb R}
\newcommand{\N}{\mathbb N}
\newcommand{\Z}{\mathbb Z}
\newcommand{\Q}{\mathbb Q}
\newcommand{\T}{\mathbb T}

\newcommand{\be}{\begin{equation}}
\newcommand{\ee}{\end{equation}}
\newcommand{\ba}{\begin{eqnarray}}
\newcommand{\ea}{\end{eqnarray}}

\def\qd{{\rule{2.3mm}{2.3mm}}}
\def\qed{{\hfill$\quad$\qd\\}}

\def\ep{{\varepsilon}}

\newtheorem{thm}{\bf Theorem}
\newtheorem{lem}[thm]{\bf Lemma}
\newtheorem{cor}[thm]{\bf Corollary}
\newtheorem{prop}[thm]{\bf Proposition}
\newtheorem{defi}[thm]{\bf Definition}
\newtheorem{rem}[thm]{\bf Remark}
\newtheorem{exa}[thm]{\bf Example}

\renewcommand{\le}{\leqslant}
\renewcommand{\ge}{\geqslant}

\date{February 3, 2010}
\begin{document}

\title{Chaotic dynamical systems associated 
with tilings of $\R^N$}

\author{
Lionel Rosier \thanks{Institut \'{E}lie Cartan,
UMR 7502 UHP/CNRS/INRIA,
B.P. 239,  54506 Vand\oe uvre-l\`{e}s-Nancy Cedex, France.
({\tt rosier@iecn.u-nancy.fr}). LR was partially
supported by the ``Agence Nationale de la Recherche'' (ANR), Project CISIFS,
grant ANR-09-BLAN-0213-02.}}

\maketitle
%
%

\begin{abstract}
In this chapter, 
we consider a class of discrete dynamical systems defined 
on the homogeneous space associated with a regular tiling of $\R^N$, whose 
most familiar example is provided by the $N-$dimensional torus $\T ^N$. 
It is proved that any dynamical system in this class
is chaotic in the sense of Devaney, and that it admits at least one positive Lyapunov exponent.
Next, a chaos-synchronization mechanism is introduced and used
for masking information in a communication setup.
\end{abstract}
\va\va
\va

\va
\section{Introduction}


Chaos synchronization has exhibited an increasing interest in
the last decade since the pioneering works 
reported in \cite{PecCar90,PecCar91}, and it has been advocated
as a powerful tool in secure communication 
\cite{Spec97a,Spec97b,KolKenChu98b,Spec00a,Spec02a}.
Chaotic systems are indeed characterized by a great sensitivity to the
initial conditions and a spreading out of the trajectories, two properties
which are very close to the Shannon requirements of confusion and diffusion
\cite{Massey92}.

There are basically two approaches when using chaotic dynamical systems for
secure communications purposes. The first one amounts to numerically
computing a great number of iterations of a discrete chaotic system,
in using e.g. the message as initial data
(see \cite{Schmitz01} and the references therein).
The second one amounts
to hiding a message in a chaotic dynamics. Only a part of the state vector
(the ``output'') is conveyed through the public channel. Next, a
synchronization mechanism is designed to retrieve the message at the
receiver part (see \cite{RMB06} and the references therein).

In both approaches, the first difficulty is to ``build'' a
chaotic system appropriate for encryption purposes.  
In this context, the corresponding chaotic signals must have
no patterning,  a broad-band power spectrum and an auto-correlation
function that quickly drops to zero.
In \cite{PecCar97}, a mean for synthesizing volume-preserving
or volume expanding maps
is provided. For such systems, there are several directions of expansion
(stretching), while the discrete trajectories are folded back into a confined
region of the phase space. Expansion can be carried out by
unstable linear mappings with at least one positive Lyapunov exponent.
Folding can be carried out with modulo functions through shift operations,
or with triangular, trigonometric functions through reflexion operations.
Fully stretching piecewise affine Markov maps have also attracted
interest because such maps are expanding in all directions and
they have uniform invariant probability densities
(see \cite{rovset98,Hasdelro96}).

Besides, we observe that the word ``chaotic'' has not the 
same meaning everywhere, and that the chaotic behavior of a 
system is often demonstrated only by numerical evidences.
The first aim of this chapter is to provide a rigorous
analysis, based on the definition given by Devaney \cite{Devaney},
of the chaotic behavior of a large class of affine dynamical systems defined 
on the homogeneous space associated with a regular tiling of $\R ^N$. 
Classical piecewise affine chaotic transformations, as the {\em tent map}, 
belong to that class. The 
dimension $N$ may be arbitrarily large in the theory developed below,  
but, for obvious reasons, most of the examples given here will be related 
to regular tilings of the plane ($N=2$).  
The study of the subclass of (time-invariant or switched) affine systems on $\T ^N$, the $N-$dimensional torus, is done in \cite{RMB04,RMB06}. 
The folding for this subclass is carried out with modulo maps, which,
from a geometric point of view, amounts to ``fold back'' 
$\R ^N$ to $[0,1)^N$ by 
means of translations by vectors in $\Z ^N$. Those translations are replaced 
here by all the isometries of some crystallographic group for an arbitrary regular 
tiling of $\R ^N$. Notice also that the fundamental domain used in 
the numerical 
implementation may be chosen with some degree of freedom. It may be a 
hypercube (as $[0,1)^N$ for $\T ^N$), or a polyhedron, or a more complicated bounded, connected set in $\R ^N$. 

For ease of implementation and duplication, a
cryptographic scheme must involve a map for which the parameters
identification is expected to be a difficult task, while computational requirements
for masking and unmasking information are not too heavy. The second aim
of this chapter is to show that all these requirements are fulfilled for
the class of dynamical systems considered here. The way of extracting the masked information is provided through an observer-based
synchronization mechanism with a finite-time
stabilization property.

Let us now describe the content of the chapter. Section 2 is devoted to the
mathematical analysis of the chaotic properties of the following
discrete dynamical system
\begin{equation}
\label{affine}
x_{k+1}=Ax_k+B \quad  (\hbox{\rm mod } G)
\end{equation}
where $A\in \Z ^{N\times N}$, $B\in \R^N$, and (mod $G$) means roughly 
that $x_{k+1}$ is the point in the fundamental domain $\cal T$ derived from
$A x_k+B$ by some transformation $g$ in the group $G$.
(\ref{affine}) may be viewed as a ``realization'' in ${\cal T}\subset \R ^N$
of an abstract dynamical system on the homogeneous space $\R ^N/G$ of classes modulo $G$.
The torus $\T ^N$ corresponds to the simplest case when 
$G$ is the group of all the translations of vectors $u\in \Z ^N$ and
the fundamental domain is ${\cal T}=[0,1)^N$. 
Note that most of the examples encountered in the literature are given only 
for the torus 
$\T ^N$ with $N=1$ and $|A|\ge 2$, or for $N=2$ and $\text{det } A=1$
(see e.g. \cite{KH}).
We give here a sufficient condition for (\ref{affine}) to be
chaotic in the sense of Devaney for any given regular tiling of $\R ^N$ ($N\ge 1$),
and we investigate the Lyapunov exponents of (\ref{affine}) and the 
equirepartition of the trajectories of (\ref{affine}).

Finally, a masking/unmasking technique based on a dynamical
embedding is proposed in Section 3.

\section{Chaotic dynamical systems and regular tilings of $\R ^N$}
\subsection{Chaotic dynamical system}
Let $(M,d)$ denote a compact metric space, and let $f:M\to M$ be a
continuous map.
The following definition of a chaotic system
is due to Devaney \cite{Devaney}.
\begin{defi}
\label{def1}
The discrete dynamical system
$$
(\Sigma )\qquad x_{k+1}=f(x_k)
$$
is said to be {\em chaotic}
if the following conditions are fulfilled:\\
{\em (C1) (Sensitive dependence on initial conditions)} There exists a number
$\varepsilon >0$ such that for any $x_0\in M$ and any $\delta >0$, there
exists a point $y_0\in M$ with $d(x_0,y_0)<\delta $  and an integer $k\ge 0$
such that $d(x_k,y_k)\ge \ep$;\\
{\em (C2) (One-sided topological transitivity)} There exists some
$x_0\in M$ with $(x_k)_{k\ge 0}$ dense in $M$;\\
{\em (C3) (Density of periodic points)} The set
$D=\{ x_0\in M; \ \exists k>0,\  x_k=x_0\}$ is dense in $M$.
\end{defi}
Recall \cite[Thm 5.9]{Walters}, \cite[Thm 1.2.2]{Vesentini}
that when $f$ is {\em onto} (i.e., $f(M)=M$), the one-sided
topological transitivity is equivalent to the condition: \\
(C$2'$)  For any pair of nonempty open sets $U,V$ in $M$,
there exists an integer $k\ge 0$ such that $f^{-k}(U)\cap V\ne \emptyset$
($\iff U\cap f^k(V)\ne \emptyset$).

\subsection{Regular tiling of $\R^N$}
\begin{figure}[hbt]
\begin{center}
\epsfig{figure=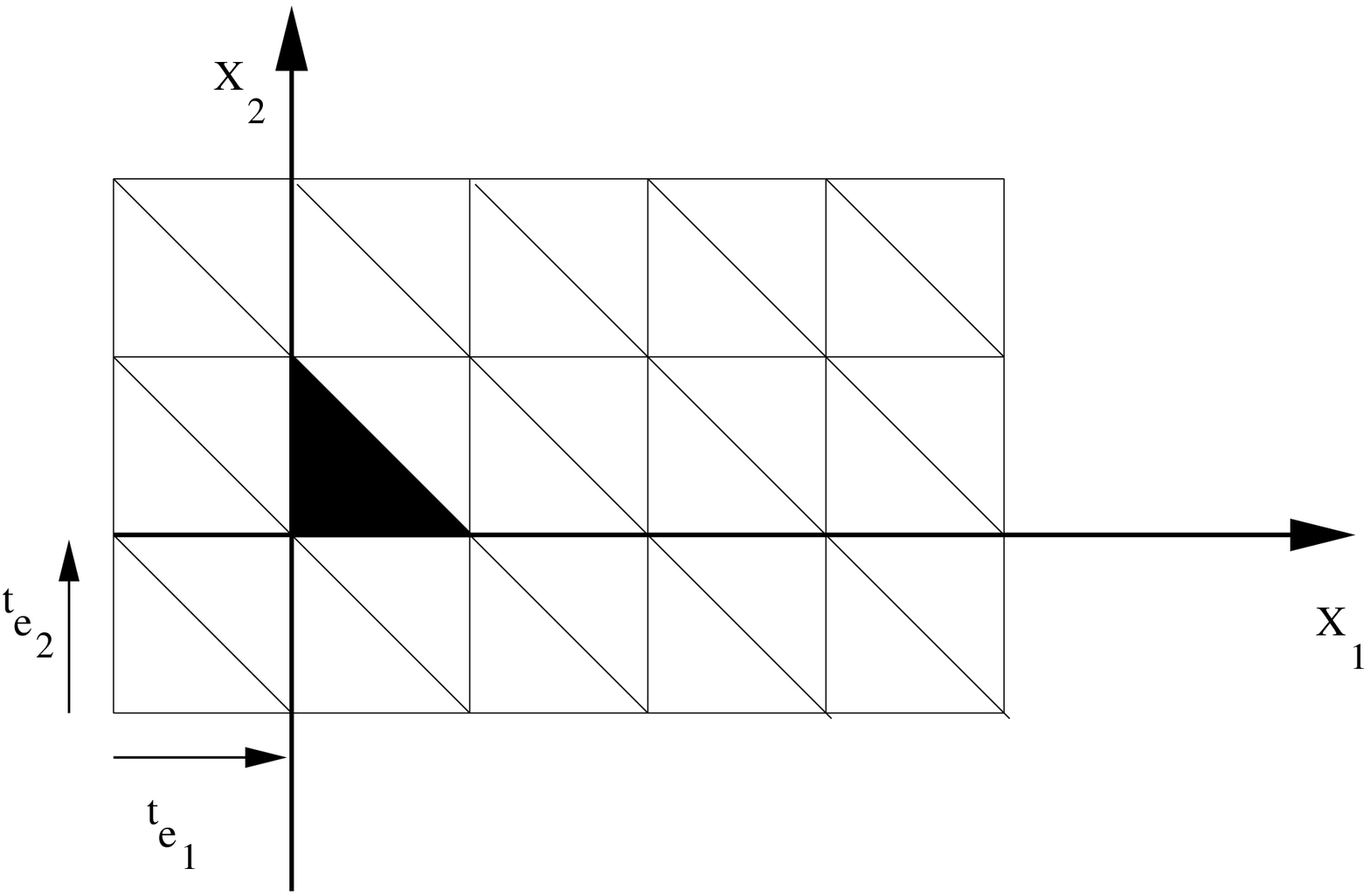,width=10cm} 
\caption{A regular tiling of $\R^2$ with a triangular fundamental tile.}
\label{triangle}
\end{center}
\end{figure}

An {\em isometry} $g$ of $\R^N$ is a map from $\R ^N$ into $\R ^N$ such that 
$||g(X)-g(Y)||=||X-Y||$ for all $X,Y\in \R ^N$. 
Let $G$ be a group of isometries of $\R ^N$ such that 
for any point $X\in \R ^N$
the orbit of $X$ under the action of $G$, namely the set
$$
G\cdot X =\{ g(X);\  g\in G \},
$$
is closed and discrete. 
Let $P\subset \R ^N$ be a compact, connected set with a nonempty interior.
Following \cite{Berger}, 
we shall say that the pair $(G,P)$ constitutes a {\em regular tiling} of $\R^N$ if
the two following conditions are fulfilled:
\begin{eqnarray}
&& \displaystyle\bigcup_{g\in G} g(P) = \R ^N \label{T1}\\
&& \forall g,h\in G\quad 
\left( g(\intP ) \cap h(\intP )\ne \emptyset \quad \Rightarrow \quad g=h\right).\label{T2}
\end{eqnarray}
Recall that $\intP$ stands for the {\em interior} of $P$, that is
$$
\intP = \{ x\in P;\ \exists \varepsilon >0,\ B(x,\varepsilon )\subset P\}.
$$ 
The set $P\subset \R ^N$ is termed a {\em fundamental tile}, and the group $G$
a {\em crystallographic group}.  
An example of a regular tiling of $\R^ 2$ with  
a triangular fundamental tile is represented in Fig. \ref{triangle}.

Note that a point
$X\in \R ^N$ may in general be obtained in several ways as the transformation
of a point in $P$ by an isometry in $G$. We introduce a set 
$\cal T$, called a {\em fundamental domain}, with $\intP \subset {\cal T}\subset P$ and such that
\begin{eqnarray}
&& \displaystyle\bigcup_{g\in G} g({\cal T}) = \R ^N \label{T10}\\
&& \forall X,X'\in {\cal T},\ \forall g\in G\quad 
\left( X'=g(X)\quad \Rightarrow \quad X'=X\right).\label{T11}
\end{eqnarray}

Introducing the equivalence relation in $\R ^N$ 
$$
X\sim Y \quad \iff \quad\exists g\in G,\quad Y=g(X),
$$
we denote by $x=\overline{X}$ the class of $X$ for $\sim$, 
i.e. $x=\{ g(X);\ g\in G\} = G\cdot X$. 
When several groups are considered at some time, we denote by $\overline{X}^G$ the class of 
$X$ modulo $G$. 
Finally, we introduce the homogeneous space of cosets $\H =(\R ^N/G) =
\{ x=\overline{X}; X\in \R ^N\}$, and define on it the following metric
$$
d(\overline{X}, \overline{Y} )=\inf_{g\in G} ||Y-g(X)||.
$$
The natural covering mapping
 $\pi : \R^N\to \H$, defined by $\pi (X)=\overline{X}$, satisfies
$$
d(\pi (X), \pi (Y)) \le ||X-Y||,
$$ 
hence it is continuous. It follows that $\H=\pi( P)$ is a compact metric space. 
On the other hand, the restriction of $\pi$ to $\cal T$ is a bijection from $\cal
T$ onto $\H$. We may therefore define the projection $\varpi:\R ^N\to {\cal T}$
by $\varpi (X)= (\pi_{\vert {\cal T}})^{-1} \pi (X)$. Note that $\varpi$ is 
in general not continuous when $\cal T$ is equipped with the topology induced from $\R ^N$, while it is continuous when $\cal T$ is endowed with the topology inherited from $\H$.

The simplest example of a regular tiling of $\R ^N$ is provided by the 
group of translations by vectors with integral coordinates 
(which is isomorphic to the lattice subgroup)
\begin{equation}
G=\{ t_u;\ u\in \Z ^N\} \sim \Z^N,
\label{Gtore}
\end{equation} 
where $t_u(X)=X+u$. 
In such a situation, a fundamental tile (resp. domain) is given by
$P=[0,1]^N$ (resp. ${\cal T}=[0,1)^N$), and the homogeneous space $\H$ is the 
standard $N-$dimensional torus $\T ^N$. 
A classification (up to isomorphism)  of the crystallographic groups 
of $\R ^N$ has been
done for a long time for $N\le 3$. There are 17 such groups in $\R ^2$, and 
230 groups in $\R ^3$, see \cite{Berger,Bur}.

\subsection{Affine transformation}
We aim to define ``simple'' chaotic dynamical systems on $M=\H$ by 
using affine transformations. Assume given a matrix $A\in \Z ^{N\times N}$
and a point $B\in \R^N$. The following hypotheses will be used at 
several places in the chapter.\\
(H1) 
$$\forall X,X'\in \R^n \quad (X\sim X' \Rightarrow 
AX+B \sim AX' +B)
$$
i.e. $X'=g(X)$ for some $g\in G$ implies $AX'+B=g'(AX+B)$ for some
$g'\in G$;\\
(H2) There exist a subgroup $G'\subset G$ of translations and a finite
collection of isometries $(g_i)_{i=1}^k$ in $G$ such that
\begin{enumerate}
\item $G$ is spanned as a group by the isometries in $G'\cup (g_i)_{i=1}^k$;
\item $G'=\{ t_u; u=\sum_{i=1}^N y_iu_i,\  y=(y_i)_{i=1}^N\in \Z^N \}$
for some basis $(u_i)_{i=1}^N$ of $\R ^N$;
\item Setting 
$     P':= \cup_{1\le i\le k} \ g_i(P)$
we have that $(G',P')$ is a regular
tiling of $\R ^N$. We denote by  ${\cal T}'$ a fundamental domain for $(G',P')$.
\end{enumerate}
(H1) is a compatibility condition needed to define a dynamical system
on $\H$. If $G$ is given by \eqref{Gtore}, then (H1) 
holds for any $A\in \Z^{N\times N}$ and any $B\in \R^N$. However, if 
\begin{equation}
\label{Greseau}
G=\{ t_u; u=\sum_{i=1}^N y_iu_i,\  y=(y_i)_{i=1}^N\in \Z^N \}
\end{equation}
for some basis $(u_i)_{i=1}^N$ of $\R ^N$, then 
(H1) holds if and only if 
\begin{equation}
\label{condAreseau}
U^{-1}AU\in \Z ^{N\times N}
\end{equation}
where $U$ is the $N\times N$ matrix with $u_i$ as $i$th column 
for $1\le i\le N$. 

(H2) allows to decompose the projection $\varpi$ onto $\cal T$ into 
a projection onto ${\cal T}'$, a fundamental domain for the 
regular tiling $(G',P')$ of $\R ^N$ involving only
translations, followed by a projection from ${\cal T}'$ onto $\cal T$. 
\begin{exa}
Let $G=<t_1,t_2,r>$ and $G'=<t_1,t_2>$, where $t_1(X)=X+(1,-1)$,
$t_2(X)=X+(1,1)$, and $r(X_1,X_2)=(-X_2,X_1)$.
Pick $k=4$ and $(g_1,g_2,g_3,g_4)=(r,r^2,r^3,id)$. 
Take as fundamental tiles $P=\{ X=(X_1,X_2); 1\le X_1\le 2,\ 
 0\le X_2\le 2-X_1\}$  (solid line) and $P'=P\cup r(P)\cup r^2(P)\cup r^3(P)$
(broken line) (see Fig. \ref{GGprime}).

\begin{figure}[hbt]
\begin{center}
\epsfig{figure=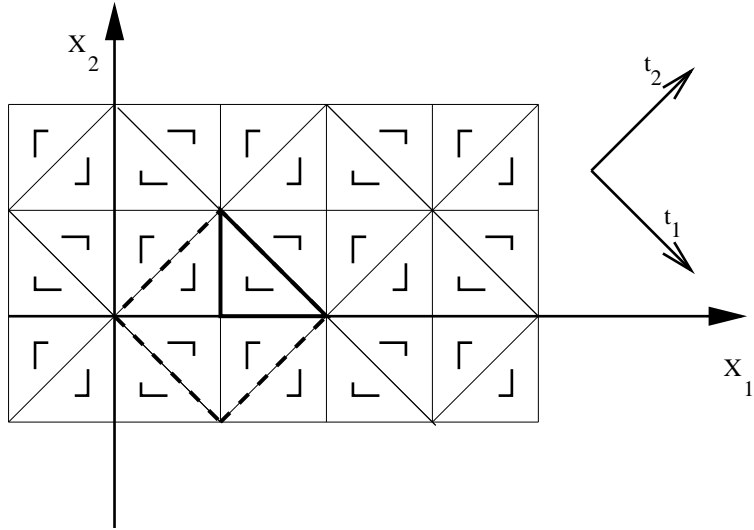,width=10cm}
\caption{A regular tiling of $\R^2$ with a triangular fundamental tile.}
\label{GGprime}
\end{center}
\end{figure}
\end{exa}

Assume that (H1) holds. Then we may define 
$$A\overline{X} +B :=\overline{AX+B}$$ 
for any $X\in \R ^N$. Thus we may consider the dynamical system
$(\Sigma _{A,B})$ on $\H$ defined by
\begin{equation}
\label{dyn}
(\Sigma _{A,B})\
\left\{
\begin{array}{rl}
x_{k+1}&=f(x_k):=A x_k+B,\\
x_0&\in \H.
\end{array}
\right.
\end{equation}

The map $f$ is called an {\em affine transformation} of $\H$.

\begin{exa}
Let $N=1$, and let $G=<t,s>$ be the group spanned by the translation
$t(X)=X+2$ and the symmetry $s(X)=2-X$. Set $P=[0,1]$. Then $(G,P)$ constitutes a regular tiling of $\R$. Note that $P$ is also a fundamental domain. Pick
$(A,B)=(2,0)\in \R ^2$. (H1) and (H2) are satisfied with $G'=\{ t_u;\ u\in 2\Z \}$, $k=2$, $g_1=s$ and $g_2=s^2=id$.
Let us write the realization of \eqref{dyn} in $P$.
Obviously, $AX\in P$ for $0\le X<1/2$, while $s(AX)=2(1-X)\in P$ for 
$1/2\le X\le 1$. Viewed in $P=[0,1]$, the dynamics reads then 
\begin{equation}
\label{tent}
x_{k+1}=h(x_k)
\end{equation}
where $h$ is the familiar tent map  (see Fig. \ref{tentfig})
\begin{figure}[hbt]
\begin{center}
\epsfig{figure=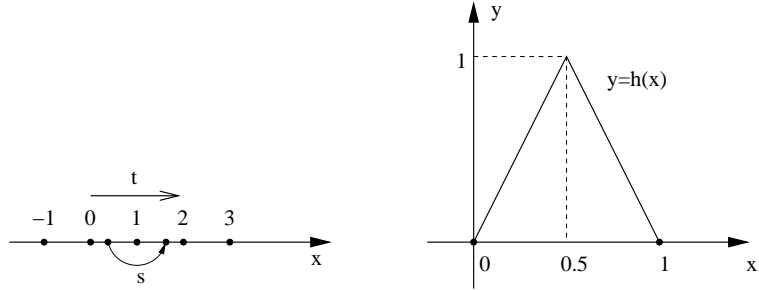,width=10cm}
\caption{A~: Action of $s$ and $t$; B~: 
the tent map}
\label{tentfig}
\end{center}
\end{figure}
$$
h(x)=\left\{ 
\begin{array}{ll}
2x \quad &\text{\rm if}\ 0\le x<\frac{1}{2},     \\
2(1-x)   &\text{\rm if}\ \frac{1}{2} \le x \le 1.
\end{array}
\right.
$$
It follows from Theorem \ref{thm2} (see below) that \eqref{tent} is chaotic on $[0,1]$.
\end{exa}
When $\H=\T^N$ and $B=0$, $f$ is nothing else than an endomorphism
of the topological group $(\T ^N,+)$, and $f$ is onto (resp., an isomorphism)
if and only if $\hbox{\rm det } A\ne 0$
(resp., $\hbox{\rm det } A=\pm 1$) (see \cite[Thm 0.15]{Walters}).
Let $\text{\rm sp}(A)$ denote the spectrum of the matrix $A$,
that is the set of the eigenvalues of $A$. A {\em root of unity}
is any complex number of the form $\lambda =\exp (2\pi it)$, with $t\in \Q$.
To see whether a dynamical system $(\Sigma _{A,B})$ is chaotic, we
need the following key result \cite[Thm 1.11]{Walters}.
\begin{prop}
\label{ergodic}
Let $f(x)=Ax+b$ ($b\in \T^N$, $A\in \Z^{N\times N}$ with
$\hbox{\rm det }A\ne 0$) be an affine transformation of $\T^N$. Then the
following conditions are equivalent:\\
\begin{tabular}{ll}
(i) &$(\Sigma _{A,b})$ is one-sided topologically transitive;\\
(ii) &(a) $A$ has no proper roots of unity (i.e., other than 1)
as eigenvalues,
and \\
&(b) $(A-I)\T ^N +\Z b$ is dense in $\T ^N$;\\
(iii) &$f$ is {\em ergodic}; that is,
$f$ is {\em measure-preserving} (i.e. for any Borel set
$E\subset \T ^N$, $m(f^{-1}(E))=$ \\
& $m(E)$, where $m$ denotes the Lebesgue measure
on $\T^N$),
and the only Borel sets $E\subset \T ^N$ for \\
&which $f^{-1}(E)=E$ satisfy $m(E)=0$ or $m(E)=1$.
\end{tabular}
\end{prop}
Notice that (ii) reduces to ``$A$ has no roots of unity as eigenvalues''
when $b=0$. Indeed, it may be seen that
$(A-I)\T ^N$ is dense in $\T ^N$ if and only if
$(A-I)$ is invertible.

\subsection{Endomorphism of $\T ^N$}
The first result in this chapter, which comes from \cite{RMB06}, provides a necessary and sufficient condition for $\Sigma _{A,0}$ to be chaotic
in $\T ^N$. 
\begin{thm}
\label{thm1}
Let $A\in \Z ^{N\times N}$. Then $(\Sigma _{A,0})$ is chaotic in 
$\T ^N$
if, and only
if, $\hbox{\rm det A}\ne 0$ and $A$ has no roots of unity as eigenvalues.
\end{thm}
{\em Proof.\ } Assume first that $(\Sigma _{A,0})$ is chaotic.
We first claim that $A$ is nonsingular.
Indeed, if $\hbox{\rm det }A=0$, then the map $f$ defined in \eqref{dyn}
is not onto \cite[Thm 0.15]{Walters}, i.e. $A\T ^N\ne \T^N$.
As $A\T ^N$ is compact (hence equal to its closure), it is not dense in
$\T ^N$, hence we cannot find some state $x_0\in \T^N$ such that
the sequence $(x_k)=(A^kx_0)$ is dense in $\T ^N$, which contradicts (C2).
Thus $\hbox{\rm det }A\ne 0$.
On the other hand, since $(\Sigma _{A,0})$ is
one-sided topologically transitive, the matrix $A$ has no roots of
unity as eigenvalues by virtue of Proposition \ref{ergodic}.

Conversely, assume that $\hbox{\rm det A}\ne 0$ and that $A$ has
no roots of unity as eigenvalues.
As (C1) is a consequence of (C2) and (C3)
(see \cite{BCDS},\cite[Thm 1.3.1]{Vesentini}), we
only have to establish the later properties. (C2) follows from Proposition
\ref{ergodic}. To prove (C3) we need to prove two lemmas.
\begin{lem}
\label{invertible}
Let $A\in \Z^{N\times N}$ be such that $\hbox{\rm det }A\ne 0$, and pick
any $p\in \N^*$ with $(p,\hbox{\rm det} A)=1$ (i.e. $p$ and
$\hbox{\rm det }A$ are relatively prime).
Then the map
$T:x\in (\Z / p\Z )^N\mapsto Ax\in (\Z/ p\Z )^N$ is invertible.
\end{lem}
{\em Proof of Lemma \ref{invertible}.\ } First, observe that the map
$T$ is well-defined. Indeed, if $X,Y\in \Z ^N$ fulfill $X-Y\in
(p\Z )^N$, then $AX-AY\in (p\Z )^N$ so that $AX$ and $AY$ belong to the same
coset in $(\Z /p\Z )^N=\Z ^N/(p\Z )^N$. As $(\Z /p\Z )^N$ is a finite set,
we only have to prove that $T$ is one-to-one. Let $X,Y\in \Z^N$ be such that
$AX=AY$ in $(\Z /p\Z )^N$ (i.e., $A(X-Y)\in (p\Z )^N$). We aim to show that
$X=Y$ in  $(\Z /p\Z )^N$ (i.e., $X-Y\in (p\Z )^N$). Set $U=X-Y$, and
pick a vector $Z\in \Z ^N$ such that $AU=pZ$. It follows that
$U=\frac{p}{\hbox{\rm det }A}{\tilde A}Z$, where
$\tilde A\in \Z^{N\times N}$ denotes the adjoint matrix of $A$ (i.e. the
transpose of the matrix formed by the cofactors).
Since $U\in \Z^N$, each component of
the vector $p{\tilde A}Z$ is divisible by
$\hbox{\rm det } A$.
Since $(p,\hbox{\rm det }A)=1$, we infer the existence of
a vector $V\in \Z ^N$ such that ${\tilde A}Z=(\hbox{\rm det }A) V$. Then
$X-Y=U=pV\in (p\Z )^N$, as desired. \qed
\begin{lem}
\label{dense}
Let $A$ and $p$ be as in Lemma \ref{invertible}, and let
$E_p:=\{\overline{0},\overline{(\frac{1}{p})}, ...,
\overline{(\frac{p-1}{p})}\}\subset \T$.  Then each point
$x\in E_p^N$ is periodic for $(\Sigma _{A,0})$. As a consequence, the set
of periodic points of $(\Sigma _{A,0})$ is dense in $\T ^N$ (i.e., (C3)
is satisfied).
\end{lem}
{\em Proof of Lemma \ref{dense}.\ }
First, observe that for any $i,j\in \{0 , ... , p-1\}$,
$i/p \equiv j/p \ (\hbox{\rm mod }1)$ if and only if  $i\equiv j\
(\hbox{\rm mod }p)$.
We infer from Lemma \ref{invertible}
that the map $\tilde T: x\in E_p^N\mapsto Ax \in E_p^N$ is well defined
and invertible. Pick any $x\in E_p^N$. As the sequence
$({\tilde T}^kx)_{k\ge 1}$
takes its values in the (finite) set $E_p^N$, there exist two numbers
$k_2>k_1\ge 1$ such that ${\tilde T}^{k_1}x={\tilde T}^{k_2}x$. $\tilde T$
being invertible, we conclude that $A^{k_2-k_1}x=x$ (i.e., $x$ is a
periodic point). Finally, the set $E=\cup \{ E_p^N;\ p\ge 1,\
(p,\hbox{\rm det }A)=1\}$ is clearly dense in $\T^N$
(take for $p$ any large prime number), and all its points
are periodic. This completes the proof of Lemma \ref{dense} and of
Theorem \ref{thm1}.\qed

For an affine transformation, we obtain a result similar to Theorem
\ref{thm1} when $1\not\in \text{\rm sp}(A)$.
\begin{cor}
\label{cor1}
Let $A\in \Z^{N\times N}$ and $b\in \T^N$. Assume that $1$ is not
an eigenvalue of $A$. Then $(\Sigma _{A,b})$ is chaotic in $\T^N$ 
if, and only if,
$\hbox{\rm det }A\ne 0$ and $A$ has no roots of unity as eigenvalues.
\end{cor}
{\em Proof.\ }
Pick any $B\in \R ^N$ with $\overline{B}=b$. As $1\not\in
\text{\rm sp} (A)$, we
may perform the change of variables
\begin{equation}
\label{cov}
x=r-\overline{(A-I)^{-1}B},
\end{equation}
which transforms (\ref{dyn}) into
\begin{equation}
\label{dynbis}
\left\{
\begin{array}{rl}
r_{k+1}&=A r_k,\\
r_0&=x_0+\overline{(A-I)^{-1}B}.
\end{array}
\right.
\end{equation}
Clearly, the conditions (C2) and (C3) are fulfilled for
$(\Sigma _{A,b})$ if, and only if,  they are fulfilled for (\ref{dynbis}).
Therefore, the result is a direct consequence of Theorem \ref{thm1}. \qed

\begin{cor}
\label{cor2}
Let $G$ be defined by \eqref{Greseau} for some basis 
$(u_i)_{i=1}^N$ of $\R ^N$. 
Let $A\in \Z^{N\times N}$ and $B\in \R ^N$. 
Assume that \eqref{condAreseau} holds
and that $1$ is not an eigenvalue of $A$. 
Then $(\Sigma _{A,B})$ is chaotic in $\H =\R ^N/G$ if, 
and only if, $\hbox{\rm det}\ A\ne 0$ and $A$ has no roots of unity as eigenvalues.
\end{cor}
{\em Proof.\ } From Corollary \ref{cor1}, we know that the dynamical system
on $\T ^N$ 
\begin{equation}
z_{k+1}={\tilde f} (z_k):=U^{-1}A U z_k + U^{-1}B 
\end{equation}
is chaotic if, and only if, ${\rm det}\ A \ne 0$ and $A$ has no roots of
unity as eigenvalues.  To prove that the dynamical system on $\H =\R ^N/G$  
\begin{equation}
x_{k+1} = f(x_k) := A x_k +B
\end{equation}
is chaotic under the same conditions, it is sufficient to prove that 
the maps $f:\H \to \H$ and $\tilde f:\T ^N\to \T ^N$ are topologically 
conjugate;
i.e., there exists a homeomorphism $h:\H \to \T ^N$ such that 
$h\circ f= \tilde f\circ h$. Define $h$ by $h(\overline{X})=\overline{Z}$
where $Z=U^{-1}X$, $\overline{X}=G\cdot X$ is the class of $X$ in $\H$
and $\overline{Z}$ is the class of $Z$ in $\T ^N$. Note first that $h$ is well
defined and continuous. Indeed, if $X'=X+UK$ with $K\in \Z ^N$, then
$Z'=U^{-1}X'=U^{-1}X+K=Z+K$, so that $h$ is well defined. On the other hand, the map $X\in \R ^N\mapsto \overline{U^{-1}X}\in \T ^N$ is clearly continuous.
Obviously, $h$ is invertible with $h^{-1}(\overline{Z})=\overline{X}$ for
$X=UZ$.  $h$ is therefore a homeomorphism from $\H$ onto $\T ^N$.  Let us check now
that $h\circ f=\tilde f\circ h$. Pick any $X\in \R ^N$. Then 
$$
h\circ f(\overline{X}^G) =h(\overline{AX+B}^{G})
=\overline{U^{-1}(AX+B)}^{\T ^N} = \tilde f (\overline{U^{-1}X}^{\T ^N}) 
=\tilde f\circ h (\overline{X}^G)
$$
and the result follows.\qed
We are in a position to state and prove the main result of this chapter. 
\begin{thm}
 \label{thm2}
Let $(G,P)$ be a regular tiling of $\R ^n$, and 
let $(A,B)\in \Z^{N\times N}\times \R ^N$ be such that both the assumptions
(H1) and (H2) are fulfilled. Assume in addition that 
$\text{\rm det}\ A\ne 0$  and that $A$ has no roots of unity as eigenvalues.  
 Then the discrete dynamical system in $\R ^N/G$
\begin{equation}
\label{A100}
x_{k+1}=Ax_k+B
\end{equation}
is chaotic.
\end{thm}
{\em Proof.}
Pick any fundamental domain $\cal T$ for $(G,P)$, and let $G'$ and ${\cal T}'$ 
be as in (H2). In addition to \eqref{A100}, we shall 
consider the discrete dynamical system
in $\R ^N/G'$
\begin{equation}
\label{A101}
z_{k+1}=Az_k+B.
\end{equation}
For any given $X_0\in \R ^N$, let $x_0=\overline{X_0}^G$  and 
$z_0=\overline{X_0}^{G'}$. Clearly, if $X\sim X'$ (mod $G'$), then 
$X\sim X'$ (mod $G$). Therefore, one can define a map $p:\R ^N/G'
\to \R ^N/G$ by $p(\overline{X}^{G'})=\overline{X}^G$. $p$ is continuous
and onto. We need two claims. \\
{\sc Claim 1.} $x_k=p(z_k)$ for all $k$. \\
Indeed, this is true for $k=0$, and if for some $k\ge 0$, $x_k=p(z_k)$  
(i.e. for some $X_k\in\R ^N$, $x_k=\overline{X_k}^G$ and $z_k=\overline{X_k}^{G'}$), 
then we have that 
$$
x_{k+1} = \overline{AX_k+B}^G = p (\overline{AX_k+B}^{G'}) = p(z_{k+1})
$$ 
which completes the proof of Claim 1.\\
{\sc Claim 2.} The image by $p$ of any dense set in $\R ^N /G'$ is a dense set
in $\R ^N/G$.\\
Let $A\subset \R ^N/{G'}$ be a given dense set. Pick any $X\in \R ^N$ and any $\varepsilon >0$.
Since $A$ is dense in $\R ^N/{G'}$, there exists $Y\in \R ^N$ such that $\overline{Y}^{G'}\in A$ and  
$$
d(\overline{X}^{G'},\overline{Y}^{G'})=\inf_{g\in G'}||Y-g(X)||<\varepsilon . 
$$
It follows that
$$
d(\overline{X}^G,\overline{Y}^G) =\inf_{g\in G}||Y-g(X)|| <\varepsilon
$$
for $G'\subset G$. Since $\overline{Y}^G = p(\overline{Y}^{G'})\in p(A)$ and the
pair $(X,\varepsilon)$ was arbitrary, 
this demonstrates that $p(A)$ is dense in $\R ^N/G$. Claim 2 is proved. 

Let us complete the proof of Theorem \ref{thm2}. To prove that \eqref{A100} is chaotic, it is sufficient (see \cite{BCDS})
to check that the conditions (C2) and (C3) are fulfilled. We know from Corollary
\ref{cor2} that \eqref{A101} is chaotic. We may therefore pick $X_0\in \R ^N$ so that, setting $z_0=\overline{X_0}^{G'}$, the sequence $\{z_k\} _{k\ge 0}$
defined by \eqref{A101} is dense in $\R ^N/G'$. By Claim 1 and Claim 2, the
sequence $\{x_k\}$ defined by \eqref{A100} and $x_0=\overline{X_0}^G$ is dense
in $\R ^N/G$; that is, (C2) is fulfilled for \eqref{A100}. On the other hand, 
the set of periodic points for \eqref{A101} is dense in $\R ^N/G'$, since 
(C3) is fulfilled for \eqref{A101}. By Claim 1, any periodic point $z_0$ for 
\eqref{A101} gives rise to a periodic point $x_0=p (z_0)$ for \eqref{A100}. 
By Claim 2, the set of periodic points for \eqref{A100} is dense in 
$\R ^N/G$; i.e., (C3) is fulfilled for \eqref{A100}. The proof of Theorem \ref{thm2}
is complete.  \qed
\begin{exa}
\begin{enumerate}
 \item Let $G=<t_{e_1},t_{2e_2},s>$ where $t_{e_1}(X)=X+(1,0)$, 
$t_{e_2}(X)=X+(0,2)$, $s(X_1,X_2)=(X_1,-X_2)$, and $P=[0,1]\times [0,1]$.
Pick $G'=<t_{e_1},t_{e_2}>$, $k=2$, $(g_1,g_2)=(s,id)$ (see Fig. \ref{sym1}).
Finally, pick
$A=\left(\begin{array}{cc}-2&0\\0&3\end{array}\right)$ and $B=(0.5, -3.2)$.
Note that $[A,S]:=AS-SA=0$, where 
$S=\left( \begin{array}{cc}1&0\\0&-1\end{array}\right)$ is the matrix 
corresponding to the symmetry $s$. 
Then (H1) and (H2) are satisfied, sp$\, (A)=\{-2,3\}$, and by Theorem \ref{thm2} 
the dynamical system \eqref{dyn} is chaotic in $\H =\R ^2/G$.

\begin{figure}[hbt]
\begin{center}
\psfig{figure=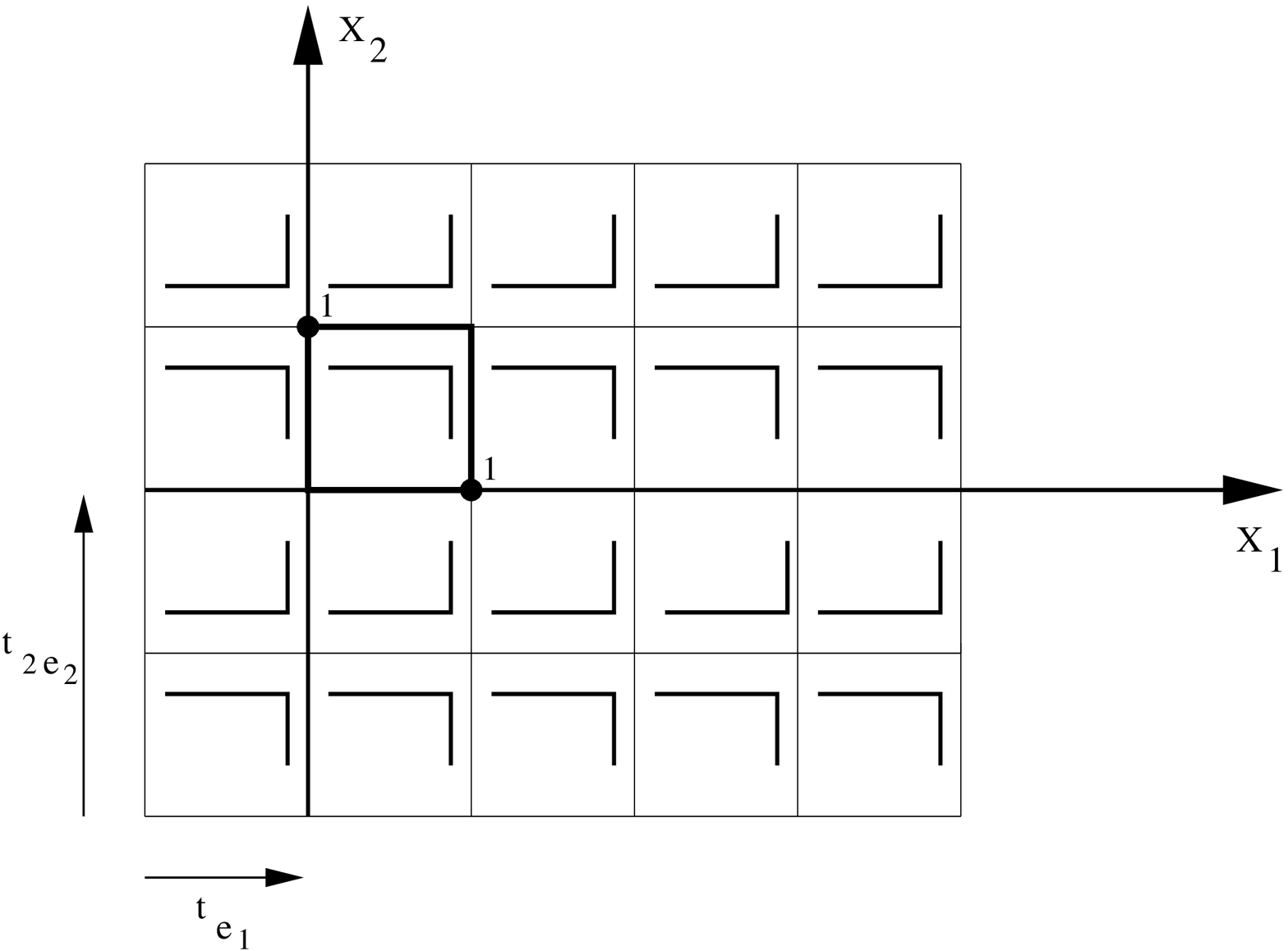,width=10cm} 
\caption{$G=<t_{e_1},t_{2e_2},s>$.}
\label{sym1}
\end{center}
\end{figure}

\item Let $G=<t_{2e_1},t_{2e_2},s_1,s_2>$ where $t_{2e_1}(X)=X+(2,0)$, 
$t_{2e_2}(X)=X+(0,2)$, $s_1(X_1,X_2)=(-X_1,X_2)$, $s_2(X_1,X_2)=(X_1,-X_2)
=-s_1(X_1,X_2)$,
and $P=[0,1]\times [0,1]$.
Pick $G'=<t_{2e_1},t_{2e_2}>$, $k=4$, $(g_1,g_2,g_3,g_4)=(s_1,s_2,s_2\circ s_1,id)$. 
(see Fig. \ref{sym2}).
Finally, pick
$A=\left(\begin{array}{cc}0&-3\\4&0\end{array}\right)$ and $B=(-0.2,1.7)$.
Note that $AS=-SA$, where $S$ is as above. 
Then (H1) and (H2) are satisfied, sp$\, (A)=\{ \pm 2i\
\sqrt{3} \}$, and by Theorem \ref{thm2} the dynamical
system \eqref{dyn} is chaotic in $\H =\R ^2/G$.

\begin{figure}[hbt]
\begin{center}
\epsfig{figure=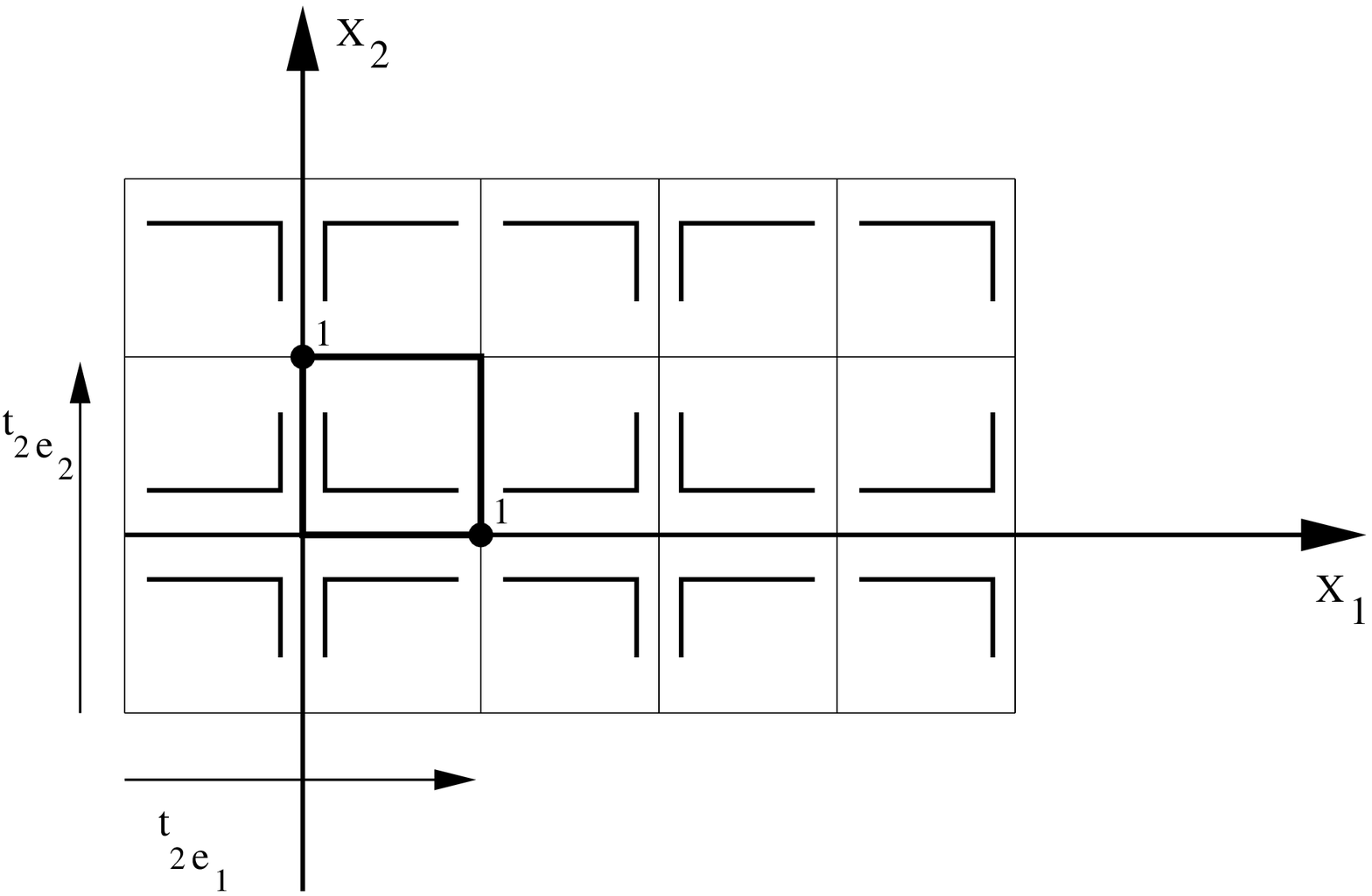,width=10cm} 
\caption{$G=<t_{2e_1},t_{2e_2},s_1,s_2>$.}
\label{sym2}
\end{center}
\end{figure}

\end{enumerate}

\end{exa}

\subsection{Lyapunov exponents}
Let $M$ denote a compact differentiable manifold endowed with a Riemann
metric $<u,v>_m$, and let $f:M\to M$ be a
map of class $C^1$. The following definition is
borrowed from \cite{Mane}.
\begin{defi}
\label{def2}
A point $x\in M$ is said to be a {\em regular point} of $f$ if there
exist numbers $\lambda _1(x)>\lambda _2(x) >\cdots > \lambda _m(x)$ and
a decomposition
$$
T_xM=E_1(x)\oplus \cdots \oplus E_m(x)
$$
of the tangent space $T_xM$ of $M$ at $x$
such that
$$\lim_{k\to +\infty} \frac{1}{k}\ln ||(D_x f^k)u||=\lambda _ j(x)$$
for all $0\ne u\in E_j(x)$ and every $1\le j\le m$.
($||v||^2:=<v,v>_x\ \forall v\in T_xM$.) The numbers
$\lambda_j(x)$ and the spaces $E_j(x)$ are termed  the
{\em Lyapunov exponents} and the {\em eigenspaces}
of $f$ at the regular point $x$.
\end{defi}
Assume now that the group $G$ is such that each isometry $g\in G$ 
has no fixed point, i.e. $g(X)\ne X$ for all $X\in \R ^N$. Then $\H=\R ^N/G$ is a smooth 
flat Riemannian manifold. Before investigating the Lyapunov exponents of an affine transformation on 
$\H$, let us give a few examples. 
\begin{exa}
 \begin{enumerate}
\item $\H =\T ^N$, and more generally, $\H =\R ^N / G$ where $G$ is as in \eqref{Greseau};
\item $\H=\R ^2 /G$ for $G=<t_{2e_1},t_{e_2},t_{e_1}\circ s>$ where $(e_1,e_2)$ is the canonical basis
of $\R ^2$ and $s(X_1,X_2)=(X_1,-X_2)$
(see Fig. \ref{klein}).
\begin{figure}[hbt]
\begin{center}
\psfig{figure=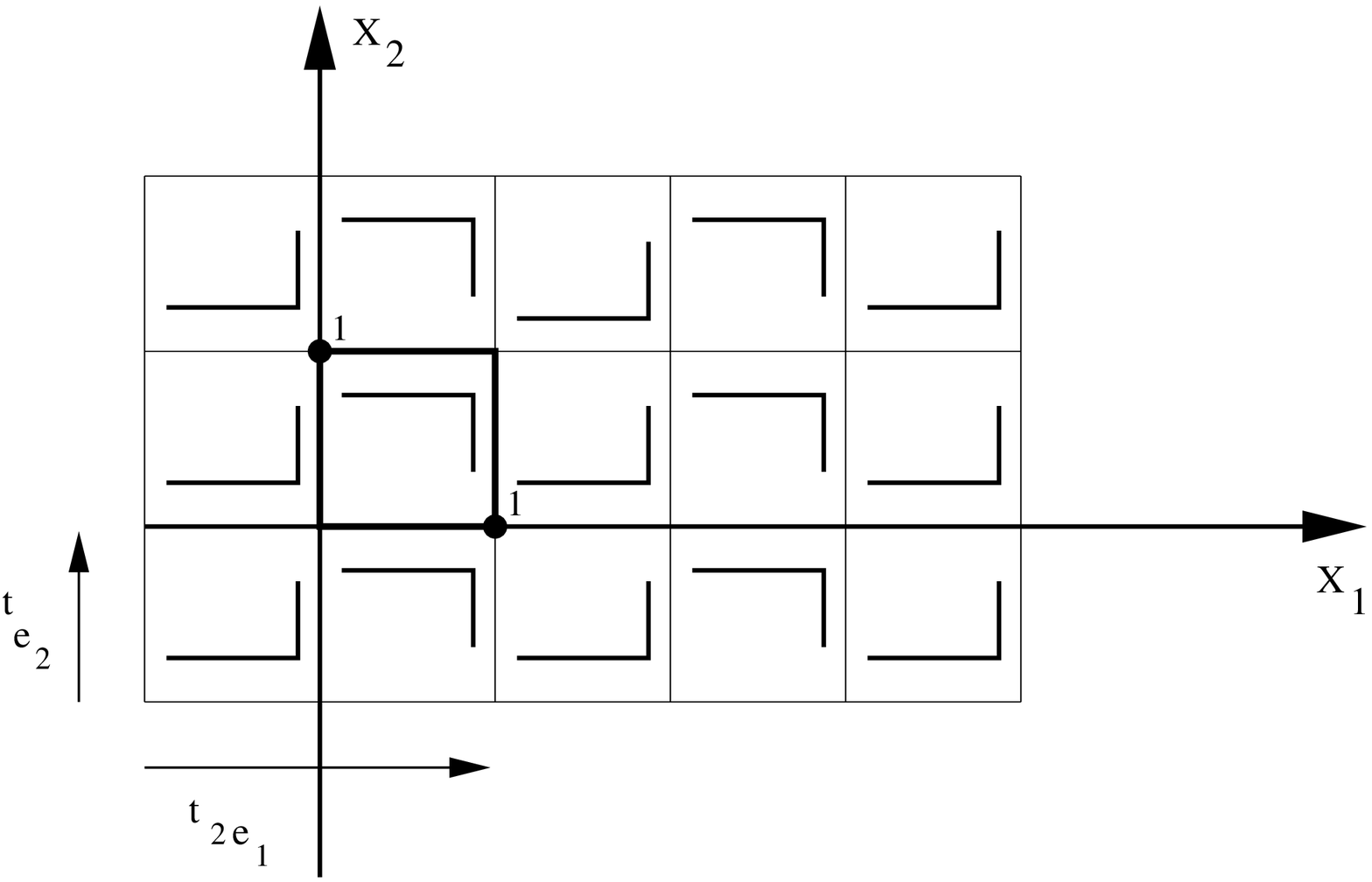,width=10cm} 
\caption{The regular tiling of $\R^2$ associated with the Klein bottle.}
\label{klein}
\end{center}
\end{figure}
$\H$ is then the {\em Klein bottle}. The torus $\T ^2$ and the 
Klein bottle $\H$ are the only smooth manifolds obtained in dimension 2. In dimension 3, there are 6 smooth
manifolds (see \cite[Section 3.5.5 p. 117]{Wolf}).  
 \end{enumerate}
\end{exa}

Consider now an affine transformation $f(\overline{X})=\overline{AX+B}$ of $\H$, 
the pair $(A,B)$ fulfilling (H1). Assume also that det $A\ne 0$. Then for any 
$k\ge 1$,
$$
f^k (\overline{X}) = \overline{A^k X + A^{k-1} B + \cdots + A B + B }.
$$
Pick a point $X\in \intP$ such that 
$$
A^k X + A^{k-1}B + \cdots + AB +  B\in \cup_{g\in G}g(\intP )
$$
(note that such a property holds for almost every $X\in \R  ^N$), and an isometry $g\in G$ such that 
$$
g(A^k X + A^{k-1}B + \cdots + AB + B )\in \intP.
$$
For $||U||$ sufficiently small, we also have that 
$$
g(A^k (X+U) + A^{k-1}B + \cdots + AB + B )\in \intP.
$$
Therefore $(D_{\overline{X}} f^k) \overline{U}=\overline{GA^kU}$, 
where $G=Dg\in \R ^{N\times N}$. Since $G$ is an orthogonal matrix, 
we have that $||\overline{GA^kU}||=||A^kU||$. 
Let $\mu _1>\mu _2>\cdots >\mu _m >0$
denote the absolute values of the eigenvalues of $A$, and let
$E_i(x)$ be the direct sum of the generalized eigenspaces
(see \cite{greub}) associated with the
eigenvalues whose absolute value is $\mu _i$, for each
$i\le m$. Then, using the Jordan decomposition of $A$,
we easily see that for any $U\in E_j\setminus \{0\}$
$$\lim_{k\to +\infty}\frac{1}{k}\ln ||A^k U|| =\ln \mu _j .$$
Observe now that if $\sigma (A)$ does not intersect the circle
$\{ z\in \C ;\ |z|=1 \}$, then $A$ has at least one eigenvalue $\lambda$ with
$|\lambda |>1$ (since the product of all the eigenvalues of $A$ is
$\hbox{\rm det } A\in \Z\setminus\{ 0\}$), hence $f$ admits at least
one {\em positive} Lyapunov exponent.
Therefore, we have proved the following
\begin{prop}
\label{prop3}
Let $(G,P)$ be a regular tiling of $\R ^N$ such that any isometry 
$g\in G$ has no fixed point. Let $(A,B)\in \R ^{N\times N} \times \R^N$ be 
such that (H1) is satisfied, $\text{\rm det}\ A\ne 0$ and each
eigenvalue $\lambda$ of $A$ satisfies $|\lambda | \ne 1$,
and let $f:\H=\R ^N/G \to \H$ be defined
by $f(x)=Ax+B$. Then almost every point $x\in\H$
is regular for $f$, with Lyapunov exponents
$\ln \mu _1 >\cdots > \ln \mu _m$, where
$\mu _1>\cdots >\mu _m$ are the absolute values of the
eigenvalues of $A$. Furthermore, $\ln \mu _1 >0$.
\end{prop}
Notice that the existence of (at least) one  positive Lyapunov exponent is often
considered as a characteristic property of a chaotic motion
\cite{WykSteeb}.
That property quantifies the sensitive dependence on initial conditions.
\subsection{Equidistribution}
In this section, $\H=\T ^N$. Let us consider a discrete dynamical system with an output
\begin{equation}
\left\{
\begin{array}{rl}
x_{k+1}&=Ax_k +B\\
y_k&=Cx_k
\end{array}
\right.
\label{output}
\end{equation}
where $x_0\in \T ^N$, $A\in \Z ^{N\times N}$, $b\in \T ^N$
and $C\in \Z ^{1\times N}$.
It should be expected that the output $y_k$
inherits the chaotic behavior of the state $x_k$.
However, Devaney's definition of a chaotic system cannot be tested
on the sequence $(y_k)$, since this sequence is not defined as
a trajectory of a dynamical system.
Rather, we may give a condition ensuring that
the sequence $(y_k)$ is equidistributed (hence dense) in $\T$ for a.e. $x_0$,
a property which may be seen as an {\em ersatz} of (C2).

If $X=(X_1, ... ,X_N),Y=(Y_1, ... ,Y_N)$ are any given points in $[0,1)^N$
and $x=\overline{X}$, $y=\overline{Y}$, then we say that $x<y$
(resp., $x\le y$) if $X_i<Y_i$ (resp., $X_i\le Y_i$)
for $i=1,...,N$. The set of points
$z\in \T^N$ such that $x\le z<y$ will be denoted by $[x,y)$.
Let $(x_k)_{k\ge 0}$ be any sequence in $\T^N$.
For any subset $E$ of $\T ^N$, let
$S_K(E)$ denote the number of points $x_k$, $0\le k\le K -1$, which
lie in $E$.
\begin{defi} \cite{KN}
We say that $(x_k)$ is {\em uniformly distributed modulo 1} (or
{\em equidistributed in $\T^N$}) if
$$
\lim_{K\to \infty} \frac{S_K([x,y))}{K} =m([x,y))=\prod_{i=1}^N(Y_i-X_i)
$$
for all intervals $[x,y)\subset \T^N$.
\end{defi}
The following result is very useful to decide whether a sequence is
equidistributed or not.
\begin{prop} ({\bf Weyl criterion} \cite{KN}, \cite{Rauzy})
The sequence $(x_k)_{k\ge 0}$ is equidistributed in $\T ^N$ if, and only if,
for every lattice point $p\in \Z ^N$, $p\ne 0$
$$
\frac{1}{K}\sum_{0\le k<K} e^{2i\pi p\cdot x_k}\to 0
\qquad \hbox{ as }K\to +\infty .
$$
\end{prop}
The next result shows that under the same assumptions as in Corollary
\ref{cor1} the sequences
$(x_k)$ and $(y_k)$ are respectively equidistributed in $\T ^N$ and $\T$
for a.e. initial state $x_0\in \T ^N$.
\begin{thm}
\label{thm3}
Let $A\in \Z ^{N\times N}$, $b\in \T ^N$ and $C\in \Z ^{1\times N}
\setminus \{ 0\}$.
Assume that $\hbox{\rm det }A\ne 0$ and that $A$ has no roots of unity as
eigenvalues (hence $\Sigma _{A,b}$ is chaotic).
Then for a.e. $x_0\in \T ^N$ the sequence $(x_k)$ (defined in (\ref{output}))
is equidistributed in $\T ^N$, and the sequence $(y_k)=(Cx_k)$ is
equidistributed in $\T$.
\end{thm}
{\em Proof:\ }
By virtue of Theorem \ref{ergodic}, the map $f(x)=Ax+b$ is ergodic
on $\T ^N$. It follows then from Birkhoff Ergodic Theorem
(see e.g. \cite[Thm 1.14]{Walters}) that for any $h\in L^1(\T ^N, dm)$ and
for a.e. $x_0\in \T^N$
$$
\frac{1}{K}\sum_{0\le k < K } h(f^k(x_0)) \to \int_{\T ^N}h(y)\, dm(y)\qquad
\hbox{\rm as } K\to +\infty .
$$
Therefore, for every lattice point $p\in \Z ^N$, $p\ne 0$, and for a.e.
$x_0\in \T ^N$
$$
\frac{1}{K}\sum_{0\le k<K}e^{2\pi i p\cdot f^k(x_0)}
\to \int_{\T ^N} e^{2\pi ip\cdot y}\, dm(y)=0\qquad
\hbox{ as } K\to +\infty .
$$
As $\Z ^N\setminus \{ 0\}$ is countable, the same property
holds for a.e. $x_0\in \T ^N$ and all $p\in \Z ^N\setminus \{ 0\}$.
Therefore, we
infer from Weyl criterion
that the sequence $(x_k)=(f^k(x_0))$ is
equidistributed for a.e. $x_0\in \T ^N$.
Pick any $x_0\in \T ^N$ such that $(x_k)$ is equidistributed, and let us
show that the output sequence $(y_k)=(Cx_k)$ is also
equidistributed provided that $C=(C_1,...,C_N)\ne (0,...,0)$.
Indeed, for any $p\in \Z \setminus \{ 0\}$
$$
\frac{1}{K}\sum_{0\le k<K} e^{2\pi i py_k}
=\frac{1}{K} \sum_{0\le k<K} e^{2\pi i (pC)\, x_k}\to 0
\qquad \hbox{ as }K\to +\infty ,
$$
hence the equidistribution of $(y_k)$ follows again by Weyl criterion.\qed
\begin{rem}
For a regular tiling $(G,P)$ of $\R ^N$, even if the sequence 
$(x_k)$ is equidistributed in $\H$, the output 
$(y_k)$ fails in general to be equidistributed in $\T$.  This is clear when
one considers a regular tiling of $\R ^2$ with the triangle 
$P=\{ X=(X_1,X_2);\quad X_1\ge 0,\ X_2\ge 0, \ X_1+X_2\le 1\}$ as fundamental tile, 
and $C=(1\quad 0)$.
\end{rem} 

\section{Synchronization and information recovering}

The aim of this section is to suggest a chaos-based encryption
scheme involving affine transformations on the homogeneous space $\H$ associated 
with some regular tiling of $\R ^ N$.
We shall provide conditions which guarantee a synchronization
with a finite-time stability of the error despite the inherent 
nonlinearity of the chaotic systems under study.

\subsection{Encryption setup}
Assume given a regular tiling $(G,P)$ of $\R ^N$ and a pair 
$(A,B)\in \R ^{N\times N}\times \R ^N$ fulfilling the assumptions of Theorem 
\ref{thm2}. For the sake of simplicity, assume further that $\R ^N/G'=\T ^N$,
so that ${\cal T}'=[0,1) ^N$. Let $\varpi : \R ^N\to {\cal T}$ and 
$\varpi ' : \R ^N\to {\cal T}'$ denote the projections on the fundamental 
domains of $(G,P)$ and $(G',P')$, respectively. Set for $k\in \N$ and $X\in \R ^N$
\begin{equation}
\label{switch}
\varpi _k(X) = 
\left\{ 
\begin{array}{ll}
\varpi '(X) & \text{ if }\ k\not\in (N+1)\N ; \\
\varpi  (X) & \text{ if }\ k\in (N+1)\N.
\end{array}
\right. 
\end{equation}

At each discrete time $k$, a symbol $m_k \in \R$ (the {\em plaintext})
of a sequence $(m_k)_{k\ge 0}$
is encrypted by a (nonlinear) encrypting
function $e$  which ``mixes'' $m_k$ and $X_k$ and produces a
{\em ciphertext} $u_k=e(X_k,m_k)$. We also assume given a decrypting
function $d$ such that  $m_k=d(X_k,u_k)$ for each $k$.
Next, the ciphertext $u_k$ is embedded in the dynamics \eqref{dyn}. We shall consider the following encryption
\begin{equation}
\label{drive1}
(\Sigma _{A,B,M,C}) \qquad \left\{
\begin{array}{rl}
X_{k+1}&= \varpi _k \{ A (X_k+M u_k) + B\} \\
Y_k&=C (X_k+M u_k)
\end{array}
\right.
\end{equation}
which corresponds to an embedding of the ciphertext
in both the dynamics and the output.
In (\ref{drive1}),
$A\in \Z ^{N\times N}$, $M\in \Z ^{N\times 1}$, and
$C\in \Z ^{1\times N}$ are given matrices, and $B\in \R ^N$.
$Y_k\in \R$ is the output conveyed to the
receiver through the channel.

$\,$From the definition of the decrypting function $d$, it is clear that
to retrieve $m_k$ at the decryption side we need to recover the pair
$(X_k,u_k)$, which in turn calls for reproducing
a chaotic sequence $(\hat{X}_k)$ synchronized   with $({X}_k)$ (i.e.,
such that $\hat{X}_k-X_k\to 0$). To this end, we propose a mechanism
based on some suitable unknown input observers, inspired from the
ones given in \cite{MilDaf03c,MilDaf04d,RMB04,RMB06}. We stress 
that the gain matrices have to be $\Z$-valued here.

For the encryption considered here, the decryption involves
the following observer-like structure
\begin{equation}
\label{response1}
(\hat{\Sigma} _{A,B,M,C}) \qquad \left\{
\begin{array}{rl}
{\hat X}_{k+1} &= \varpi_k \{ A {\hat X}_k+L (Y_k-{\hat Y}_k) + B \} \\
{\hat Y}_k&=C{\hat X}_k
\end{array}
\right.
\end{equation}
where $L\in \Z ^{N\times 1}$, ${\hat X}_k\in \R ^N$ and ${\hat Y}_k\in \R$
(${\hat X}_0$ being an arbitrary point in $\R ^N$). 
Let $\overline{X}$ denote the class of $X$ modulo $G'$, i.e. in $\T ^N$.
Set $e_k=\overline{X_k}-\overline{{\hat X}_k}$ for all $k\ge 0$.
Noticing that for all $X\in \R ^N$ 
$$
\overline{\varpi _k (X)} = \overline{\varpi '(X)} = \overline{X}\qquad 
\text{ for } 1\le k\le N,
$$
we obtain by subtracting (\ref{response1})
from (\ref{drive1}) that the error dynamics reads
\begin{equation}
\label{errordyn}
e_{k+1}=(A -L C )e_k+\overline{(A-L C ) M u_k}, \qquad 
1\le k\le N.
\end{equation}

Before proceeding to the design of the observers, we give a few definitions
and a preliminary result.
\subsection{Definitions and preliminary results}
\begin{defi}
A pair $(A ^\flat, C^\flat )$ is said to be in a
{\em companion canonical form} if it takes the form
\begin{equation}
\label{companion}
A^\flat =
\left(
\begin{array}{lcccc}
- \alpha ^{N-1}        & 1         & 0         & \cdots &  0 \\
-  \alpha ^{N-2}   & 0         & 1         & \cdots &  0 \\
\vdots             & \vdots    & \vdots    & \ddots &  \vdots \\
- \alpha ^1         & 0         & 0         & \cdots &  1      \\
-\alpha ^0         & 0         & 0         & \cdots &  0
\end{array}
\right), \qquad
C^\flat =
\left(
\begin{array}{ccccc}
1 &
0 &
\cdots &
0 &
0
\end{array}
\right)\cdot
\end{equation}
\end{defi}
It is well known that the characteristic polynomial of
$A^\flat$ reads
$\chi_{A^\flat}(\lambda )
=\lambda ^N+\alpha ^{N-1}\lambda ^{N-1}+\cdots +\alpha ^1\lambda +\alpha ^0$.
\begin{defi}
Two pairs $(A,C)$ and $(A^\flat, C^\flat )$ in
$\Z ^{N\times N}\times \Z ^{1\times N}$ are said to be {\em similar over
$\Z$} if there exists a matrix $T\in \Z ^{N\times N}$ with
$\hbox{\rm det }T=\pm 1$ (hence $T^{-1}\in \Z ^{N\times N}$ too) such that
\begin{equation*}
A=T^{-1}A^\flat T,\quad C=C^\flat T.
\end{equation*}
\end{defi}

The following result provides a sufficient condition for an observable pair
$(A,C)$ to admit a $\Z$-valued gain matrix $L$ such that $A-LC$ is Hurwitz.
\begin{prop}
\label{prop1}
Let $A\in \Z^{N\times N}$ and $C\in \Z ^{1\times N}$ be two
 matrices such that $(A,C)$ is
similar over $\Z$ to a pair
$(A^\flat ,C^\flat )\in \Z^{N\times N}\times \Z ^{1\times N}$ in a companion canonical form. Let us denote by $(-\alpha^{N-1}~\cdots~-\alpha^0)'$ the first
column of $A^\flat$.
Then there exists a unique
matrix $L\in \Z ^{N\times 1}$ such that the matrix $A-LC$ is
Hurwitz (i.e., $\text{sp}(A-LC) \subset \{ z\in \C ;\ |z|<1\}$), namely
$L=T^{-1}L^\flat$ with $L^\flat=(-\alpha^{N-1}~\cdots~-\alpha^0)'$.
Furthermore, $(A-LC)^N=0$.
\end{prop}
{\em Proof.\ } Write $A=T^{-1}A^\flat T$, $C=C^\flat T$, with
$(A^\flat ,C^\flat )$ as in (\ref{companion}) and $T\in \Z^{N\times N}$
with $\hbox{\rm det }T=\pm 1$. For any given matrix
$L\in  \Z^{N \times 1}$,  we define the matrix
$L^\flat =(l^{N-1} \cdots \ l^0 )'$ by $L^\flat =TL$. Then,
$A-LC=T^{-1}(A^\flat -L^\flat C^\flat)T$ with
$$
A^\flat -L^\flat C^\flat =\left(
\begin{array}{lcccc}
- \alpha ^{N-1}-l^{N-1}         & 1         & 0         & \cdots &  0 \\
- \alpha ^{N-2}-l^{N-2}    & 0         & 1         & \cdots &  0 \\
\vdots    & \vdots       & \vdots    & \ddots    & \vdots \\
 -\alpha ^1-l^1       & 0         & 0         & \cdots & 1        \\
- \alpha ^0-l^0            & 0         & 0         & \cdots & 0
\end{array}
\right) .
$$
Its characteristic polynomial reads
$$\chi _{A^\flat  -L^\flat C^\flat}(\lambda )=
\lambda ^{N}+(\alpha ^{N-1}+l^{N-1})\lambda ^{N-1}+\cdots
+(\alpha ^1+l^1)\lambda
+(\alpha ^0+l^0).$$
If $L$ is such that $A-LC$ is Hurwitz, then $A^\flat -L^\flat C^\flat
=T(A-LC)T^{-1}$ is Hurwitz too, hence
we may write
$\chi _{A-LC}(\lambda )=\chi _{A^\flat -L^\flat C^\flat}(\lambda )=
\lambda ^p \chi(\lambda )$,
where $p\in \{ 0,...,N \}$
and $\chi \in \Z [\lambda ]$ has its roots $\lambda _1,...,\lambda _{N-p}$ in
the set $\{ z\in \C;\ 0<|z|<1\}$. Assume that $p<N$, and denote by $q$
the constant coefficient of $\chi$. Then $q\ne 0$ (since $\chi (0)\ne 0$),
and $|q|=\prod _{i=1}^{N-p}|\lambda _i| <1$, which is impossible, since
$q\in \Z$. Therefore $p=N$ and $l^j=-\alpha ^j$ for any $j\in \{0,...,N-1\}$
(hence $L^\flat$ and $L$ are unique). On the other hand
\begin{equation}
\label{eq100t2}
A^\flat -L^\flat C^\flat =\left(
\begin{array}{lcccc}
0        & 1         & 0         & \cdots &  0 \\
0    & 0         & 1         & \cdots &  0 \\
\vdots    & \vdots       & \vdots    & \ddots    & \vdots \\
0       & 0         & 0         & \cdots & 1        \\
0            & 0         & 0         & \cdots & 0
\end{array}
\right) \cdot
\end{equation}
For this choice of $L$, $\chi _{A-LC}(\lambda )
=\lambda ^{N}$ and $(A-LC)^{N}=0$. \qed
It should be emphasized that the above argument shows that
a $\Z$-valued matrix ${\cal N}$
is Hurwitz if and only if it is nilpotent. In other words, the
system $\nu _{k+1}={\cal N} \nu _k$ is {\em asymptotically stable}
if and only if it is {\em finite-time stable}.

We are now in a position to state the second main result of this chapter. 
\begin{thm}
\label{thm4}
Let $(G,P)$ be a regular tiling of $\R ^N$, and let 
$(A,B)\in \Z ^{N\times N}\times \R ^N$ be such that (H1) and (H2) 
are fulfilled with $\R ^N/G'=\T ^N$. Assume given $C\in \Z ^{1\times N}$ such that 
$(A,C)$ is similar over $\Z$ to a pair $(A^\flat, C^\flat )$ in a companion canonical form. Then one can pick two matrices $L\in \Z ^{N\times 1}$ and
$M\in \Z ^{N\times 1}$ so that 
$(A-LC)M=0$ and $CM=1$. Furthermore 
$$
X_k=\hat X_k \text { and } u_k=Y_k-\hat Y_k \qquad \forall k\ge N+1.
$$
\end{thm}
{\em Proof.} Let $T$, $A^\flat$, $C^\flat$, $L$ and 
$L^\flat$ be as in the proof of Proposition \ref{prop1}. Set 
$M^\flat =(1\ 0 \ \cdots 0 )'$ and $M=T^{-1}M^\flat $. Then 
$(A-LC)M=T^{-1}(A^\flat -L^\flat C^\flat )T\cdot T^{-1}M^\flat =0$  
by \eqref{eq100t2}, and $CM=C^\flat T\cdot T^{-1}M^\flat =1$. On the other hand, it follows from \eqref{errordyn} and the choice of $M$ that 
$$
e_{k+1}=(A-LC)e_k \qquad \forall k\in \{ 1, ..., N \}
$$
hence $e_{N+1}=(A-LC)^Ne_1=0$. Since $X_{N+1}$ and $\hat X_{N+1}$ belong
to ${\cal T}'$ by construction, we have that $\hat X_{N+1}=X_{N+1}$. 
To complete the proof, it is sufficient to prove the following\\ 
{\sc Claim.} For any $k\ge 0$, $\hat X_k=X_k$ implies $\hat X_{k+1}=X_{k+1}$.\\
Indeed, using the fact that $(A-LC)M=0$ and $\hat X_k=X_k$ we obtain that
\begin{eqnarray*}
\hat X_{k+1} 
&=& \varpi _k (A\hat X_k + L C (X_k + M u_k -\hat X_k ) +B) \\
&=& \varpi _k ( A X_k + A M u_k +B) \\
&=& X_{k+1}. 
\end{eqnarray*}
This completes the proof of Theorem \ref{thm4}. \qed  

\begin{rem}

\begin{enumerate}
 \item The projection $\varpi _k(x)$ allows to switch between 
the dynamics \eqref{A100} and \eqref{A101} in $\R /G$ and $\R /G'$, 
respectively. For a dynamics in $\T ^N$ only ($G'=G$), one can replace
$\varpi _k(x)$ by $\varpi '(x)$ (the projection onto $[0,1)^N$).
\item The result in Theorem \ref{thm4} remains true if we take 
$\varpi_k(x)=\varpi '(x)$ for $k\le N$ and $\varpi_k(x)=\varpi (x)$ for $k\ge N+1$.
However, the definition of $\varpi_k(x)$ in \eqref{switch} guarantees that 
a finite time synchronization occurs even if the output $Y_k$ is
not transmitted at some times.  Such a  property may be useful for the secured transmission of video sequences.
\item The output $Y_k=C(X_k+Mu_k)$ may be replaced by $\tilde Y_k=h(Y_k)$, 
where $h:\R \to \R$ is a nonlinear invertible map. This renders the analysis
of the dynamics of $Y_k$ much more complicated. 
\item In practice, when $\H =\T ^N$, the matrices $A,C,L$ and $M$ may be constructed in the following way. Pick any matrix $\hat T =[\hat T_{i,j}]\in \Z ^{N\times N}$ with $\hat T_{i,j}=0$ for $i>j$ and $\hat T_{i,i}=1$ for all $i$. 
We set $T=\hat T'\ \hat T$. Note that $\text{det} \ \hat T=
\text{det} \ T=1$. Next, we pick a pair $(A^\flat, C^\flat )$ in a companion canonical form so that the roots of $\chi _{A^\flat}$ do not belong to 
the set $\{ 0\}\cup \{ z\in \C;\ |z|=1\}$. Then $A,C,L$ and $M$ are defined by 
$$
A=T^{-1}A^\flat T,\quad C=C^\flat T, \quad L=T^{-1}A^\flat (C^\flat )',\quad
\text{ and } \quad M=T^{-1}(C^\flat )'.
$$
\end{enumerate}
\end{rem}

\subsection{Numerical simulations} 
\begin{figure}[hbt]
\begin{center}
\epsfig{figure=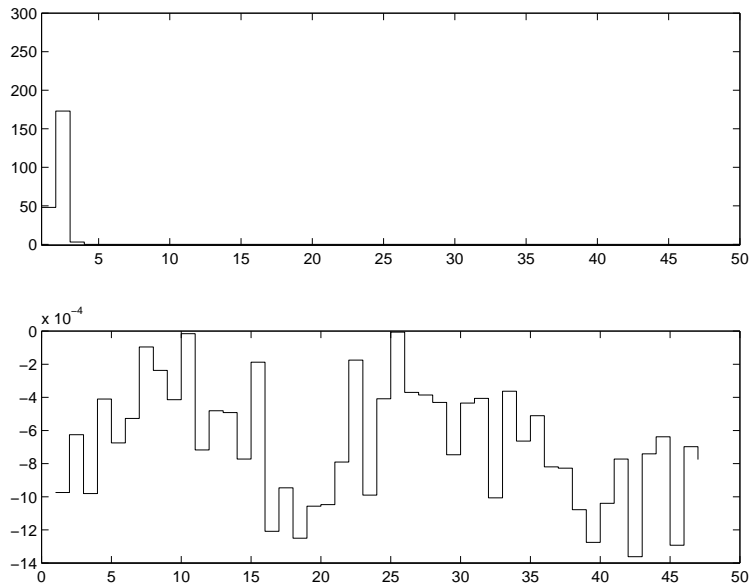,width=10cm} 
\caption{A~: error on the recovered information $u_k-\hat{u}_k$; B~: 
state reconstruction error $X_k-\hat{X}_k$}
\label{fig_1}
\end{center}
\end{figure}

This section is borrowed from \cite{RMBpreprint}. 
Assume $\H=\T ^3$ and consider the dynamical system $(\Sigma _{A,b,M,C})$ with
$$
A=
\left(
\begin{array}{ccc}
-19 & 26 & 7\\
-51 & 65 & 17\\
152 & -184 & -47
\end{array}
\right),~~C=(6~-5~-1),~~b=0.
$$
$(\Sigma _{A,b})$ is chaotic by virtue of Theorem~\ref{thm1}, 
since $\hbox{\rm det } A=3$ (hence det $A\ne 0$) and the eigenvalues of $A$ are 
-3, -0.4142, 2.4142 ($A$ has no roots 
of unity as eigenvalues). The pair $(A,C)$ is similar over 
$\Z$ to the pair $(A^\flat,C^\flat )$ in companion 
canonical form, where 
$$
A^{\flat}=
\left(
\begin{array}{ccc}
-1 & 1 & 0\\
7 & 0 & 1\\
3 & 0 & 0
\end{array}
\right),\ 
C^{\flat}=(1\ 0\ 0)
\  \hbox{\rm and }
~~T=
\left(
\begin{array}{ccc}
6 & -5 & -1\\
-5 & 10 & 3\\
-1 & 3 & 1
\end{array}
\right) .
$$
According to Proposition~\ref{prop1}, the unique matrix 
$L\in \Z ^{N\times 1}$ such that $A-LC$ is Hurwitz is 
$L=T^{-1}L^{\flat}$, with $L^{\flat}=(-1\ 7\ 3)^T$. We obtain
$L=(-2\ -6\ 19)^T$. The corresponding matrix 
$M\in \Z ^{3\times 1}$ such that $(A-LC)M=0$ and $CM=1$ is 
$M=(1\ 2\ -5)^T$.\\
The information to be masked is a flow corresponding to integers 
ranging from 0 to 255. The data are scaled to give an 
input $u_k$ ranging from 0 to $1$, and are embedded 
into the chaotic dynamics of $(\Sigma_{A,b,M,C})$. From a practical point of 
view, the transmitted signal $y_k$ cannot be coded with an infinite 
accuracy and so it has to be truncated for throughput purpose. The observer 
($\hat{\Sigma}_{A,b,M,C}$) is used in order to recover the information.   
Numerical experiments bring out that the number of digits of the conveyed 
output can actually be limited without giving rise to recovering errors. 
The results  
reported in Fig.~\ref{fig_1} show a perfect recovering for a number 
of digits of $y_k$ equal to 4 (this is the minimum number required for 
perfect retrieving). 
The recovering error reaches zero after 3 steps, a fact which is consistent 
with above theoretical results on finite time synchronization ($N=3$).
The figure highlights the fact that even though the state reconstruction may not be
perfect (residual errors due to truncations), a perfect information reconstruction
is nevertheless achieved.
\begin{rem}
Actually, for {\em any} system $\Sigma _{A,B,M,C}$, the numerical computations
can be performed {\em in an exact way}, i.e. without rounding errors, provided that
the number of digits is sufficiently large. 
\end{rem}


\subsection{Concluding remarks}
The {\em message-embedding} masking technique studied here
does not originate from the conventional cryptography 
(see \cite{Menal96} for a good survey). Nevertheless, it
seems to be highly related to some popular encryption schemes, the so-called
{\em stream ciphers} \cite{Millamigo05a}.  
Therefore, it is desirable that the proposed scheme
be robust against both statistical and algebraic attacks. On one hand, the
robustness against statistical attacks follows from the chaotic behavior of the
output. On the other hand, the security against algebraic attacks rests on
the difficulty to identify the parameters of the 
system. The identification of the parameters 
is here a hard task for two reasons:
\begin{enumerate}
\item
The particular {\em structure} of the encryption system
$(\Sigma _{A,B,M,C})$,
that is the {\em dimension} of the matrix $A$ and the tiling of the space
used, is assumed to be unknown;
\item The ciphertext $u_k$ actually results from a mixing between the
plaintext $m_k$ and the state $X_k$ ($u_k=e(X_k,m_k)$).
This generally results in a {\em nonlinear} dynamics
$(\Sigma _{A,B,M,C})$, rendering the parameters hardly identifiable
\cite{LjungGlad94}.
\end{enumerate}
A real-time
implementation has already been carried out on an experimental platform
involving a secured multimedia communication. (For details about the
platform, see e.g. \cite{Milal03}). 



\end{document}